\documentclass[11pt,a4paper]{article}
\usepackage{amssymb}          

\newtheorem{thm}{Theorem}[section]
\newtheorem{prop}[thm]{Proposition}

\newcommand{\eqref}[1]{(\ref{#1})}
\newcommand{\bea}{\begin{eqnarray}}
\newcommand{\ena}{\end{eqnarray}}
\newcommand{\be}{\begin{eqnarray*}}
\newcommand{\en}{\end{eqnarray*}}
\newcommand{\lb}[1]{\label{#1}}

\newcommand{\mZ}{{\mathcal Z}}
\newcommand{\F}{{\mathcal F}}

\newcommand{\id}{{\rm id}}
\newcommand{\tr}{{\rm tr}}

\newcommand{\ve}{{\varepsilon}}

\newcommand{\z}{{\zeta}}

\newcommand{\ket}[1]{{| #1 \rangle}}      

\newcommand{\Z}{{\mathbb Z}}
\newcommand{\C}{{\mathbb C}}

\newcommand{\g}{\mathfrak{g}}
\newcommand{\slh}{\widehat{\mathfrak{sl}}_2}
\newcommand{\slN}{\widehat{\mathfrak{sl}}_N}
\newcommand{\slf}{\mathfrak{sl}_2}

\begin{document}
\font\csc=cmcsc10 scaled\magstep1

\vspace*{21mm}
\begin{center}
{\Large\bf 
On Lepowsky-Wilson's ${\mathcal Z}$-algebra
} 
\end{center}
\vspace{11mm}
\setcounter{footnote}{0}\renewcommand{\thefootnote}{\arabic{footnote}}
\begin{center}
{\csc
Yuji Hara\footnote{Graduate School of Arts and Sciences, 
The University of Tokyo, Komaba, Tokyo 153-8914, Japan.\\
e-mail:snowy@gokutan.c.u-tokyo.ac.jp},
Michio Jimbo\footnote{Graduate School of Mathematical Sciences,
The University of Tokyo, Komaba, Tokyo 153-8914, Japan.\\
e-mail:jimbomic@ms.u-tokyo.ac.jp}, 
Hitoshi Konno\footnote{
Department of Mathematics, Faculty of Integrated Arts and Sciences,
Hiroshima University, \\ Higashi-Hiroshima 739-8521, Japan.
e-mail:konno@mis.hiroshima-u.ac.jp},	\\
Satoru Odake\footnote{
Department of Physics, Faculty of Science, Shinshu University, 
Matsumoto 390-8621, Japan.\\
e-mail:odake@azusa.shinshu-u.ac.jp}
and 
Jun'ichi Shiraishi\footnote{Graduate School of Mathematical Sciences,
The University of Tokyo, Komaba, Tokyo 153-8914, Japan.\\
e-mail:shiraish@ms.u-tokyo.ac.jp}

}
\end{center}

\begin{abstract}
We show that 
the deformed Virasoro algebra specializes in a certain limit to 
Lepowsky-Wilson's ${\mathcal Z}$-algebra. 
This leads to a free field realization of the affine Lie algebra 
${\widehat{\mathfrak{sl}}_2}$ which respects the principal gradation. 
We discuss some features of this bosonization including the screening 
current and vertex operators. 
\end{abstract}
\newpage
\section{Introduction}

In 1978,  Lepowsky and Wilson introduced the idea of 
free field representation in the theory of affine Lie algebras. 
In the first paper \cite{LW}, a representation 
of the affine Lie algebra $\slh$ was found for the special case of level one. 
Soon later this construction was extended to more general cases 
by introducing an associative algebra 
called the $\mZ$-algebra \cite{z-alg},  
and a connection between the Rogers-Ramanujan identities 
and the affine Lie algebras was uncovered. 
This study, however, was restricted to the case of the level $k$  
being a non-negative integer, 
because their idea of representing the $\mZ$-algebra 
was to use $k$-copies of Heisenberg algebras. 

We now know that the Wakimoto construction \cite{Wak86} 
affords a way of realizing representations of an arbitrary level $k$. 
The Wakimoto representation plays an important role,  
both as a powerful computational tool in mathematical physics,  
as well as for theoretical studies in pure representation theory. 
There is, however, one basic difference between the Wakimoto 
and Lepowsky-Wilson's original constructions. 
While the former is based on the homogeneous Heisenberg subalgebra,   
the latter is based on the principal Heisenberg subalgebra.  
Though strange as it may seem, 
we have not seen in the literature a construction 
which respects the principal gradation and 
works for an arbitrary level.

In this note we revisit this problem. 
We observe here a simple but peculiar fact that the $\mZ$-algebra 
arises as a certain specialization of the deformed Virasoro algebra 
\cite{FR95,qVir}.
The representation of $\slh$ mentioned above is obtained as 
an immediate consequence of this observation. 
We also touch upon the screening currents, vertex operators and a 
connection to the elliptic Knizhnik-Zamolodchikov (KZ) equation 
studied by Etingof \cite{Et94}. 


\section{Free field realization of $\slh$}

\subsection{$\slh$ in the principal picture} 

In order to fix the notation, let us recall the realization of 
 $\slh$ in the principal picture. 
Let $e,f,h$ be the standard generators of $\g=\slf$ with
$[h,e]=2e$, $[h,f]=-2f$, $[e,f]=h$.    
The invariant bilinear form on $\slf$ is chosen as $(e,f)=1,(h,h)=2$. 
Let $\g_0=\C h$, $\g_1=\C e\oplus\C f$. 
The affine Lie algebra $\slh$ is realized as a vector space 
\be
\slh=\g_0\otimes\C[t^2,t^{-2}]\oplus
\g_1\otimes t\C[t^2,t^{-2}]
\oplus \mathbb{C} c\oplus \mathbb{C} \rho,
\en
endowed with the Lie bracket 
\be
&&[a\otimes t^m ,b\otimes t^n ]=[a,b]\otimes t^{m+n}+
 {1\over 2}(a,b) mc\;\delta_{m+n,0},\\
&&c\,:\mbox{\,central},\quad 
[\rho,a\otimes t^m]=m\,a\otimes t^m,
\en
where $a,b\in\g$. 
We set 
\be
&&
\beta_n=(e+f)\otimes t^n \quad (\mbox{$n$ odd}), 
\qquad
x_n=
\cases{
h\otimes t^n & (\mbox{$n$ even}),\cr
(-e+f)\otimes t^n & (\mbox{$n$ odd}).\cr
}
\en
In terms of the generating series (currents) 
\be
\beta(\zeta)=\sum_{n:{\rm odd}}	
\beta_n \zeta^{-n},\qquad
x(\zeta)=\sum_{n \in \mathbb{Z}} x_{n} \zeta^{-n}, 
\en
the commutation relations read as 
\bea
&&[\beta(\zeta),\beta(\xi)]={c\over 2}
\left((D\delta)\left({\xi\over \zeta}\right)
-(D\delta)\left(-{\xi\over \zeta}\right) \right),
\lb{eq1}
\\
&&[\beta(\zeta),x(\xi)]=
\left(\delta\left({\xi\over \zeta}\right)
-\delta\left(-{\xi\over \zeta}\right) \right)x(\xi),
\lb{eq2}\\
&&[x(\zeta),x(\xi)]=
-2 \delta\left(-{\xi\over \zeta}\right)\beta(\xi)
+c  (D\delta)\left(-{\xi\over \zeta}\right).
\lb{eq3}
\ena
Here $\delta(\zeta)=\sum_{m\in \mathbb{Z}} \zeta^m$ and 
$D=D_\zeta$ stands for $\zeta{d\over d\zeta}$. 

In the sequel we fix a complex number $k\neq -2$ (the level),   
and focus attention to representations of $\slh$ 
on which the central element $c$ acts as $k$ times the identity. 
Our aim is to find a free field realization of the relations 
\eqref{eq1}-\eqref{eq3}. 
For that purpose let us introduce three kinds of bosonic fields 
\be
&&\phi_1(\zeta)=-\sum_{n:{\rm even}\atop n\neq 0} 
{\phi_{1,n}\over n}\zeta^{-n}+\phi_{1,0}\log\zeta+Q,
\\
&&
\phi_i(\zeta)=-\sum_{n:{\rm odd}}	
{\phi_{i,n}\over n}\zeta^{-n}
\qquad (i=0,2). 
\en
The fields $\phi_0(\zeta)$, $\phi_2(\zeta)$ are odd,  
$\phi_i(-\zeta)=-\phi_i(\zeta)$, 
while the derivative of $\phi_1(\zeta)$ is even, 
$(D \phi_1)(-\zeta)=(D \phi_1)(\zeta)$.  
We set the commutation relations for their Fourier modes as 
\be
&&[\phi_{1,m},\phi_{1,n}]=4(k+2)m \delta_{m+n,0},\qquad
[\phi_{1,0},Q]=4(k+2),\\
&&
[\phi_{0,m},\phi_{0,n}]=4km \delta_{m+n,0},\\
&&
[\phi_{2,m},\phi_{2,n}]=-4km \delta_{m+n,0}, 
\en
all other commutators being $0$. 
For $j\in\C$, we denote by $\F_{j,k}$ the Fock space for three bosons  
\be
&&
\F_{j,k}=\mathbb{C}[\,\phi_{0,-n},\phi_{2,-n}~~(n=1,3,5,\cdots),
~~\phi_{1,-n}~~(n=2,4,6,\cdots)\,]\ket{j,k}
\en
generated on the Fock vacuum  $\ket{j,k}$,
\be
\phi_{i,n}\ket{j,k}=0\quad (n> 0, i=0,1,2),
\quad \phi_{1,0}\ket{j,k}=2j\,\ket{j,k},
\quad e^{\frac{Q}{2(k+2)}}\ket{j,k}=\ket{j+1,k}.
\en
We denote by $d:\F_{j,k}\rightarrow\F_{j,k}$ the grading operator
\be
[d,\phi_{i,n}]=-n\,\phi_{i,n},\quad 
[d,Q]=\phi_{1,0},
\qquad d\ket{j,k}=\frac{2j^2+k}{4(k+2)}\ket{j,k}.
\en
We adopt the conventional `normal ordering rule' 
and the `normal ordering symbol' $:\cdots:$ for our bosonic fields.
For example,
\begin{eqnarray*}
:e^{\phi_1(\zeta)}:
=e^{\sum_{n\geq 1}{\phi_{1,-2n}\over 2n}\zeta^{2n}}
e^{-\sum_{n\geq 1}{\phi_{1,2n}\over 2n}\zeta^{-2n}}
e^Q \zeta^{\phi_{1,0}}.
\end{eqnarray*}

We can now state the main result of this note. 
\begin{prop}
Let $j,k$ be complex numbers with $k\neq 0,-2$. 
Then the following gives a level $k$ representation 
of $\slh$ on the Fock space $\F_{j,k}$: 
\bea
&&\beta(\zeta) ={1\over 2}D\phi_0(\zeta),
\lb{eq7}
\\
&&
x(\zeta)={1\over 2}:
\left( D\phi_1(\zeta)+D\phi_2(\zeta) \right)
e^{{\phi_2(\zeta)\over k}+{\phi_0(\zeta)\over k}}:,
\lb{eq8}\\
&&c=k,\quad \rho=-d. 
\lb{eq9}
\ena
\end{prop}

This representation is highest weight in the sense that 
\be
&&\beta_n \ket{j,k}=0,~x_n \ket{j,k}=0\quad {\rm for}\; n>0,\\
&&x_0 \ket{j,k}=j\ket{j,k}. 
\en
The highest weight is $\frac{k}{2}(\Lambda_1+\Lambda_0)+
j(\Lambda_1-\Lambda_0)$ where $\Lambda_0,\Lambda_1$ are 
the fundamental weights of $\slh$.  
The character of this representation 
(counted according to the principal gradation)
\be
\mbox{tr}_{\F_{j,k}}(q^{-\rho})
=q^{\frac{2j^2+k}{4(k+2)}}
\frac{1}{(q;q)_\infty (q;q^2)_\infty},
\en
is the same as that of the Verma module of $\slh$. 
Here $(z;p)_\infty=\prod_{n=0}^\infty(1-p^nz)$.  
In particular, the above representation is irreducible for generic 
values of $j,k$. 

\subsection{Connection with the deformed Virasoro algebra}

Because of the commutation relation \eqref{eq1},  
the current $\beta(\zeta)$ can tautologically be identified with 
the bosonic field $(1/2)D\phi_0(\zeta)$.  
As Lepowsky and Wilson have shown, 
the other current $x(\zeta)$ can be realized as 
\bea
&&
x(\zeta)= z(\zeta):e^{\phi_0(\zeta)\over k} :, 
\lb{eq10}
\ena
provided $z(\zeta)=\sum_{n\in\Z} z_n\zeta^{-n}$ commutes with $\phi_0(\zeta)$ 
and satisfies the relation of the $\mZ$-algebra
\bea
\left(\zeta_1-\zeta_2\over \zeta_1+\zeta_2\right)^{2/k}
z(\zeta_1)z(\zeta_2)
=
\left(\zeta_2-\zeta_1\over \zeta_2+\zeta_1\right)^{2/k}
z(\zeta_2)z(\zeta_1)+ k (D\delta)\left(-{\zeta_2\over \zeta_1}\right).
\lb{eq5}
\ena
As we explain below,  
this algebra is related with the deformed Virasoro algebra (DVA). 

The DVA is an associative algebra generated by $T_n$ ($n\in\Z$) 
(see \cite{qVir}). 
In terms of $T(\zeta)=\sum_{n\in \mathbb{Z}} T_n \zeta^{-n}$ 
the defining relations read 
\bea
f(\zeta_2/\zeta_1)T(\zeta_1)T(\zeta_2)-T(\zeta_2)T(\zeta_1)f(\zeta_1/\zeta_2)
= -\frac{(1-q)(1-t^{-1})}{1-p}\left[
       \delta \Bigl(\frac{p\zeta_2}{\zeta_1}\Bigr)-
       \delta \Bigl(\frac{p^{-1}\zeta_2}{\zeta_1}\Bigr)\right],
\lb{eq4}
\ena
where  $q$ and $t$ are parameters, $p=q/t$, and 
\be
f(\zeta)
=\exp \left\{\sum_{n\geq 1}\frac{1}{n}\frac{(1-q^n)(1-t^{-n})}{1+p^n}
 \zeta^n \right\}.
\en
Now set 
\begin{eqnarray*}
&&q=e^h,\qquad t=-q^{k+2\over 2}, 
\end{eqnarray*}
and consider the limit $h\rightarrow 0$. 
Suppose that the expansion 
\bea
T(\zeta)= 0+h T^{(1)}(\zeta) +h^2T^{(2)}(\zeta) +\cdots
\lb{eq6}
\ena
takes place. 
Under this assumption we find that, at the second order in $h$, 
the relation \eqref{eq4} reduces to the $\mZ$-algebra relation 
\eqref{eq5} with the identification $\sqrt{-1}T^{(1)}(\zeta)=z(\zeta)$. 
We have verified the expansion \eqref{eq6} 
using the known bosonization for DVA \cite{qVir}. 
This leads to the formula 
\be
z(\zeta)=\frac{1}{2}:(D\phi_1(\zeta)+D\phi_2(\zeta))
e^{\frac{\phi_2(\zeta)}{k}}:. 
\en
The formula \eqref{eq8} 
for $x(\zeta)$ follows from this and \eqref{eq10}.
In fact this bosonization formula for the $\mZ$-algebra 
was first obtained by a guesswork 
based on the one for the ordinary Wakimoto realization 
which respects the homogeneous gradation \cite{Frau}. 
The relationship between the $\mZ$-algebra 
and the deformed Virasoro algebra was noticed only afterwards.

Thus, somewhat unexpectedly,   
one can view the DVA with the parameters $q$ and $t=-q^{k+2\over 2}$ 
as a quantum deformation of Lepowsky-Wilson's $\mZ$-algebra at level $k$.  

The DVA admits two kinds of ``screening operators'' $S_{\pm}(\xi)$ 
commuting with $T(\zeta)$ up to a total difference. 
In the next section we shall show that one of them
in the limit, $S(\xi)$, becomes 
the screening operator for the present bosonization of $\slh$. 
The limit of the other one $\eta(\xi)$, after a modification by zero-mode,  
plays a role of the operator ``$B$'' which appears 
in the construction of the elliptic KZ equation \cite{Et94}.
\medskip

\noindent{\it Remark.}\quad 
It seems likely that a similar construction persists in the case of  
the deformed $W_N$ algebra \cite{FeFr95,AKOS96}. 
Let $\omega$ be a primitive $N$-th root of unity, and  
set $t=\omega q^{\frac{k+N}{N}}$. 
Suppose that in the limit $q\rightarrow 1$   
we have the expansion of the $W$-currents 
\bea
W_i(\zeta)= 0+h W_i^{(1)}(\zeta) +h^2 W_i^{(2)}(\zeta) +\cdots.
\lb{eq20}
\ena
We have checked for $N=3$ (and partially for all $N$) 
that the $W_i^{(1)}$ ($i=1,\cdots,N-1$) then 
satisfy the relations of the $\mathcal{Z}$-algebra for $\slN$. 
However, for $N\ge 3$
we have not been able to show \eqref{eq20}, 
which seems to hold not in the free field realization 
but only at the level of correlation functions. 

\section{Vertex operators and KZ equations}
\subsection{Screening current $S(\xi)$}
Our representation of $\slh$ on the 
Fock space $\F_{j,k}$ may become reducible for some specific values 
of the parameters $j,k$. 
In the usual Wakimoto construction 
the object which controls this phenomenom is the screening current. 
Let us consider its analog. 
Define $S(\zeta):\F_{j,k}\to\F_{j-2,k}$ by
\begin{eqnarray*}
&&S(\zeta)={1 \over  2}\zeta^{2 \over  k+2}:
D\phi_2(\zeta)
e^{-{\phi_1(\zeta)\over k+2}}:.
\end{eqnarray*}
It enjoys the expected properties 
\be
&&
[\beta(\zeta),S(\xi)]=0,\\
&&
[x(\zeta),S(\xi)]\\
&&\quad=
{k+2\over 2}
D_\xi
\left(-
\delta\left(\xi\over \zeta\right) \xi^{\frac{2}{k+2}}
:e^{-{\phi_1(\xi)\over k+2}+{\phi_2(\xi)\over k}+
{\phi_0(\xi)\over k}}:
+
\delta\left(-{\xi\over \zeta}\right) \xi^{\frac{2}{k+2}}
:e^{-{\phi_1(\xi)\over k+2}-{\phi_2(\xi)\over k}-
{\phi_0(\xi)\over k}}:
\right).
\en
These equations imply that the screening charge 
\be
&&
Q^m=
\oint\cdots 
\oint \frac{d\xi_1}{\xi_1}\cdots \frac{d\xi_m}{\xi_m}
S(\xi_1)\cdots S(\xi_m):
\F_{j,k}\longrightarrow \F_{j-2m,k},
\en
commutes with the action of $\beta_n,x_n$, if closed contours
for all the $\xi_i$'s exist.

\subsection{Vertex Operators}
Let us consider the vertex operator (VO) 
associated with the two-dimensional representation. 
Set $V=\C u_+\oplus \C u_-$, and let 
$V'(\z)=V\otimes\mathbb{C}[\z,\z^{-1}]$ be the $\slh$-module given by 
\be 
&&
\beta_n (u_\pm\otimes \zeta^m) =\mp u_\pm \otimes \zeta^{m+n},
\qquad
x_{n} (u_\pm\otimes \zeta^m )
=\cases{
- u_\mp \otimes \zeta^{m+n} & ($n$ even), \cr
\mp u_\mp  \otimes \zeta^{m+n}& ($n$ odd). \cr
}
\en
We define $\Phi(\zeta):\,\F_{j,k}\rightarrow \F_{j+1,k}\otimes V'(\z)$ 
by 
\be
&&
\Phi(\zeta)v=\Phi_+(\zeta)v\otimes u_++
\Phi_-(\zeta)v\otimes u_-,\\
&&\Phi_\pm(\zeta)=\zeta^{{1\over 2(k+2)}}
:e^{{\phi_1(\zeta)\over 2(k+2)}\pm
{\phi_2(\zeta)\over 2k}\pm
{\phi_0(\zeta)\over 2k}}:.
\en
Then we have the intertwining property 
\be
(x\otimes\id+\id\otimes x) \Phi(\zeta)
=\Phi(\zeta)x\quad ({}^\forall x\in \slh). 
\lb{eqn:def-VO}
\en
\medskip

\noindent {\it Remark.}\quad 
$V'(\zeta)$ contains a proper submodule 
 $W=\mbox{span}\{u_+\otimes \zeta^m-u_-\otimes (-\zeta)^m
\mid m\in\Z\}$.  
Its quotient $V(\zeta)=V'(\zeta)/W$ is isomorphic to the irreducible 
evaluation module of $\slh$ associated with $V$. 
The above VO naturally gives rise to the intertwiner
 $\F_{j,k}\rightarrow \F_{j+1,k}\otimes V(\z)$. 

\subsection{Elliptic Knizhnik-Zamolodchikov equation}

As usual, the highest-to-highest matrix elements of VO's satisfy 
the KZ equation. 
Let us consider the elliptic KZ equation studied by Etingof \cite{Et94}. 

Let $M_{j,k}$ denote the Verma module over $\slh$ with highest weight
 $\frac{k}{2}(\Lambda_1+\Lambda_0)+j(\Lambda_1-\Lambda_0)$
and highest weight vector $v_{j,k}$.  
Denote by $\Psi(\zeta):M_{j,k}\rightarrow M_{j',k}\otimes V(\zeta)$ the 
intertwining operator. 
Let $B:M_{j,k}\to M_{-j,k}$ denote the linear isomorphism characterized by 
\be
&&
Bv_{j,k}=v_{-j,k}, \\
&&
B\beta(\zeta)=\beta(\zeta)B,\qquad Bx(\zeta)=-x(\zeta)B.
\en
It was shown in \cite{Et94} that the trace function
\be
&&F(\zeta_1,\cdots,\zeta_n)
=\tr_{M_{j,k}}\left(\Psi(\z_1)\cdots\Psi(\z_n)Bq^{-\rho}\right)\\
\en
satisfies the following elliptic KZ equation 
\begin{eqnarray*}
&&
(k+2)D_{\zeta_i}F(\zeta_1,\cdots,\zeta_n)
=
\sum_{j(\ne i)}r^{ij}(\zeta_i/\zeta_j)
F(\zeta_1,\cdots,\zeta_n), 
\en
where $r(z)$ denotes the elliptic classical $r$ matrix 
\be
&&
r(z)=
\left(\sum_{l\in\mathbb{Z}}{q^l z \over (q^l z)^2-1}\right)
(e+f)\otimes(e+f)
-
\left(\sum_{l\in\mathbb{Z}}{(-q)^l z \over (q^lz)^2-1}\right)
(e-f)\otimes (e-f)
\\
&&
\qquad
+\left(\sum_{l\in\Z}(-)^l\frac{(q^lz)^2+1}{(q^lz)^2-1}\right)
\frac{1}{2}h\otimes h,
\en
and $r^{ij}(z)$ signifies $r(z)$ acting nontrivially 
on the $(i,j)$-th components.

Let us look for a bosonic realization of the operator $B$. 
In the usual Wakimoto construction, 
we have the `$\xi$-$\eta$ system' after bosonizing the `$\beta$-$\gamma$
ghost system'.  
An analog of $\eta$ in the present case is 
\be
&&
\eta(\zeta)
=\sum_{n\in\mathbb{Z}}\eta_n\z^{-n}
=\zeta^{k+2\over 2}
:e^{{\phi_1(\zeta)\over 2}+ {\phi_2(\zeta)\over 2}} :,
\en
which satisfies 
\begin{eqnarray*}
&&
[\beta(\zeta),\eta(\xi)]=0,\\[2mm]
&&
[x(\zeta),\eta(\xi)]_+=
2\xi \partial_\xi
\left(
\xi^{k+2 \over 2} \delta\left(\frac{\zeta}{\xi}\right)
:e^{{ \phi_1(\xi)\over 2}+ {(k+2)\phi_2(\xi)\over 2k}+\frac{\phi_0(\xi)}{k}}:
\right),
\end{eqnarray*}
where $[A,B]_+=AB+BA$. 
{}From the above, we find that the zero-th Fourier mode 
$\eta_0:\F_{j,k}\to\F_{j+k+2,k}$ satisfies 
\be
\eta_0\beta(\zeta)=\beta(\zeta)\eta_0,\qquad\eta_0x(\zeta)=-x(\zeta)\eta_0.
\en
Therefore, if $j+k+2=-j+2m$ for some $m\in\Z_{\ge 0}$, 
the combination of $\eta_0$ and the screening operator $Q^m$ implements $B$. 

In particular, when $m=0$ and $2j\in \Z_{>0}$, 
the screening operators do not appear and 
the traces can be expressed without using integrals. 
For example, when $(j,k)=(1/2,-3)$ and $(j,k)=(1,-4)$ we find
\be
&&\tr_{\F_{1/2,-3}}
\left(\Phi_\pm(\zeta_1)\,\eta_0\,q^{-\rho}\right)
=q^{\frac {5}{8}}(q^2;q^2)_\infty^{-\frac {3}{2}},
\\
&&
\tr_{\F_{1,-4}}\left(\Phi_{\ve_1}(\z_1)\Phi_{\ve_2}(\z_2)\,\eta_0\,
q^{-\rho}\right)
\\[2mm]
&&
\qquad =q^{\frac 14}(q^2;q^2)_\infty^{-\frac 94}
\z^{1/4}\Theta_{q^2}(\z^2)^{-\frac 14}\Theta_{q^2}(-\ve_1\ve_2q\z)\\[2mm]
&&
\qquad =q^{\frac{1}{2}}
\left\{\frac{\sqrt{-1}\varepsilon_1\varepsilon_2}{8\pi^3}
\wp'\!\left(\frac{\ln(-\ve_1\ve_2q\z)}{2\pi\sqrt{-1}}\Bigm|1,2\tau\right) 
\right\}^{-\frac 14},
\en
where $\z=\z_2/\z_1$, $q=\exp(2\pi\sqrt{-1}\tau)$, and
$\wp(z|\omega,\omega')$ denotes the Weierstrass elliptic 
function with fundamental periods $\omega,\omega'$. 
These formulas have been discussed in \cite{Et94} 
(there are minor errors in Section 5, equation (5.2) of \cite{Et94}). 

\end{document}